# Spectral sets for locally bounded operators[1]

Stoian Sorin Mirel

**Abstract**

We introduce the spectral radius $r_\mathcal{P}(T)$ for a quotient bounded operator on a locally convex space X. Similarly to the case of bounded operator on a Banach space we prove that the Neumann series $\sum_{n=0}^{\infty} \frac{T^n}{\lambda^{n+1}}$ converges to $R(\lambda,T)$, whenever $|\lambda| > r_\mathcal{P}(T)$, and $|\sigma(Q_\mathcal{P}, T)| = r_\mathcal{P}(T)$. For locally bounded operators on sequentially complete locally convex space we have
$$|\sigma(T)| = |\sigma_{lb}(T)| = |\sigma(Q,T)| = r_{lb}(T),$$
and the spectral set $\sigma(T)$ is compact.

## 1. Introduction

The spectral theory for a linear on Banach space X is well developed and we have useful tools for use this theory. For example, the spectral radius of such operator T is defined by the Gelfand formula $r(T) = \lim_{n \to \infty} \sqrt[n]{\|T^n\|}$ and $|\sigma(Q,T)| = r(T)$. Further it is known that the resolvent $R(\lambda,T)$ is given by the Newmann series $\sum_{n=0}^{\infty} \frac{T^n}{\lambda^{n+1}}$, whenever $|\lambda| > r(T)$.

If we want to generalize this theory on locally convex space X one major difficulty is that is not clear which class of operators we can use, because there are several non-equivalent ways of defining bounded operators on X. The concept of bounded element of a locally convex algebra was introduced by Allan [1]. An element is said to be bounded if some scalar multiple of it generates a bounded semi group.

**Definition 1.1.** Let X be a locally convex algebra. The radius of boundness of an element $x \in X$ is the number
$$\beta(x) = \inf\{\alpha > 0 \mid \text{the set } \{(\alpha x)^n\}_{n \geq 1} \text{ is bounded}\}$$

Through this paper all locally convex spaces will be assumed Hausdorff, over complex field **C**, and all operators will be linear. If X and Y are topological vector spaces we denote by L(X,Y) ($\mathcal{L}(X,Y)$) the algebra of linear operators (continuous operators) from X to Y

Any family $\mathcal{P}$ of seminorms who generate the topology of locally convex space X (in the sense that the topology of X is the coarsest whit respect to which all seminorms of $\mathcal{P}$ are continuous) will be called a calibration on X. A calibration is characterized by the property that for every seminorms $p \in \mathcal{P}$ and every constant $\varepsilon > 0$ the sets
$$S(p, \varepsilon) = \{x \in X \mid p(x) < \varepsilon\},$$

---

[1] This work is supported by MEdC-ANCS grant ET 65/2005

constitute a neighborhoods sub-base at 0. A calibration on X will be principal if it is directed. The set of calibration for X is denoted by $\mathbf{C}(X)$ and the set of principal calibration is denoted by $\mathbf{C}_0(X)$.

Any family of seminorms on a linear space is partially ordered by relation „$\leq$", where
$$p \leq q \iff p(x) \leq q(x), (\forall) \, x \in X.$$
A family of seminorms is preordered by relation „$\prec$", where
$$p \prec q \iff \text{there exists some } r > 0 \text{ such that } p(x) \leq r\, q(x), \text{ for all } x \in X.$$
If $p \prec q$ and $q \prec p$, we write $p \approx q$.

**Definition 1.2.** Two families $\mathcal{P}_1$ and $\mathcal{P}_2$ of seminorms on a linear space are called Q-equivalent (denoted $\mathcal{P}_1 \approx \mathcal{P}_2$) provided:
a) for each $p_1 \in \mathcal{P}_1$ there exists $p_2 \in \mathcal{P}_2$ such that $p_1 \approx p_2$;
b) For each $p_2 \in \mathcal{P}_2$ there exists $p_1 \in \mathcal{P}_1$ such that $p_2 \approx p_1$.

It is obvious that two Q-equivalent and separating families of seminorms on a linear space generate the same locally convex topology.

Similary to the norm of an operator on a normed space V. Troistsky [10] define the mixed operator seminorm of an operator between locally convex spaces. If $(X,\mathcal{P})$, $(Y,\mathcal{Q})$ are locally convex spaces, then for each $p, q \in \mathcal{P}$ the application $m_{pq}: L(X,Y) \to \mathbf{R} \cup \{\infty\}$, defined by
$$m_{pq}(T) = \sup_{p(x) \neq 0} \frac{q(Tx)}{p(x)},$$
is called the mixed operator seminorm of T associated with $p$ and $q$. When X=Y and $p = q$ we use notation $\hat{p} = m_{pp}$.

**Lemma 1.3.** ( [10]) If $(X,\mathcal{P})$, $(Y,\mathcal{Q})$ are locally convex spaces and $T \in L(X,Y)$, then
1) $m_{pq}(T) = \sup_{p(x)=1} q(Tx) = \sup_{p(x) \leq 1} q(Tx)$, $(\forall) \, p \in \mathcal{P}$, $(\forall) \, q \in \mathcal{Q}$;
2) $q(Tx) \leq m_{pq}(T) p(x)$, $(\forall) \, x \in X$, whenever $m_{pq}(T) < \infty$.

**Corollary 1.4.** If $(X,\mathcal{P})$, $(Y,\mathcal{Q})$ are locally convex spaces and $T \in L(X,Y)$, then
$$m_{pq}(T) = \inf\{M > 0 \mid q(Tx) \leq M\, p(x), (\forall) x \in X\},$$
whenever $m_{pq}(T) < \infty$.

**Proof.** If $p, q \in \mathcal{P}$ then from previous lemma we have
$$q(Tx) \leq m_{pq}(T) p(x), (\forall)\, x \in X.$$
If $M > 0$ such that
$$q(Tx) \leq M\, p(x), (\forall)\, x \in X.$$
then using (1) we obtain
$$m_{pq}(T) = \sup_{p(x)=1} q(Tx) \leq M.$$

**Definition 1.5.** An operator T on a locally convex space X is quotient bounded with respect to a calibration $\mathcal{P} \in \mathbf{C}(X)$ if for every seminorm $p \in \mathcal{P}$ there exists some $c_p > 0$ such that
$$p(Tx) \leq c_p\, p(x), (\forall)\, x \in X.$$

The class of quotient bounded operators with respect to a calibration $\mathcal{P} \in \mathbf{C}(X)$ is denoted by $Q_{\mathcal{P}}(X)$.

**Lemma 1.6.** If $X$ is a locally convex space and $\mathcal{P} \in \mathbf{C}(X)$, then for every $p \in \mathcal{P}$ the application $\hat{p}: Q_{\mathcal{P}}(X) \to \mathbf{R}$ defined by
$$\hat{p}(T) = \inf\{ r > 0 \mid p(Tx) \leq r\, p(x), (\forall)\, x \in X \},$$
is a submultiplicative seminorm on $Q_{\mathcal{P}}(X)$, satisfying $\hat{p}(I) = 1$.

We denote by $\hat{\mathcal{P}}$ the family $\{ \hat{p} \mid p \in \mathcal{P} \}$.

**Proposition 1.7.** ([7]) Let $X$ is a locally convex space and $\mathcal{P} \in \mathbf{C}(X)$.
1) $Q_{\mathcal{P}}(X)$ is a unital subalgebra of the algebra of continuous linear operators on $X$;
2) $Q_{\mathcal{P}}(X)$ is a unital local multiplicative convex algebra (l.m.c. -algebra) with respect to the topology determined by $\hat{\mathcal{P}}$.
3) If $\mathcal{P}' \in \mathbf{C}(X)$ such that $\mathcal{P} \approx \mathcal{P}'$, then $Q_{\mathcal{P}'}(X) = Q_{\mathcal{P}}(X)$ and $\hat{\mathcal{P}} = \hat{\mathcal{P}}'$;
4) The topology generated by $\hat{\mathcal{P}}$ on $Q_{\mathcal{P}}(X)$ is finer than the topology of uniform convergence on bounded subsets of $X$.

**Lemma 1.8.** If $X$ is a sequentially complete convex space, then $Q_{\mathcal{P}}(X)$ is a sequentially complete m-convex algebra for all $\mathcal{P} \in \mathbf{C}(X)$.

**Proof.** Let $\mathcal{P} \in \mathbf{C}(X)$ and $(T_n)_n \subset Q_{\mathcal{P}}(X)$ be a Cauchy sequence. Then, for each $\varepsilon > 0$ and each $\hat{p} \in \hat{\mathcal{P}}$ there exists some index $n_{p,\varepsilon} \in \mathbf{N}$ such that
$$|\hat{p}(T_n) - \hat{p}(T_m)| \leq \hat{p}(T_n - T_m) < \varepsilon,\ (\forall)\, n, m \geq n_{p,\varepsilon}. \tag{1}$$

From previous relation it follows that $(\hat{p}(T_n))_n$ is convergent sequence of real numbers, for each $\hat{p} \in \hat{\mathcal{P}}$. If $x \in X$, then
$$p(T_n x - T_m x) \leq \hat{p}(T_n - T_m)\, p(x),\ (\forall)\, p \in \mathcal{P}, \tag{2}$$
so $(T_n(x))_n \subset X$ is a Cauchy sequence. But, since $X$ is sequentially complete and separated, there exists an unique element $y \in X$ such that
$$\lim_{n \to \infty} T_n x = y.$$

Therefore, the operator $T: X \to X$ defined by
$$T(x) = \lim_{n \to \infty} T_n x,\ (\forall)\, x \in X,$$
is well defined. It is obvious that $T$ is linear operator. Using the continuity of seminorms $\hat{p} \in \hat{\mathcal{P}}$ we have
$$p(Tx) = p(\lim_{n \to \infty} T_n x) = \lim_{n \to \infty} p(T_n x) \leq \lim_{n \to \infty} \hat{p}(T_n)\, p(x) = c_p\, p(x),$$
for all $x \in X$ and for each $p \in \mathcal{P}$ (where $c_p = \lim_{n \to \infty} \hat{p}(T_n)$). This implies that $T \in Q_{\mathcal{P}}(X)$.

Now we prove that $T_n \to T$ in $Q_{\mathcal{P}}(X)$. From relations (1) and (2) it follows that for each $\varepsilon > 0$

and each $\hat{p} \in \hat{\mathcal{P}}$ there exists $n_{p,\varepsilon} \in \mathbf{N}$ such that
$$p(T_n x - T_m x) \leq \varepsilon\, p(x),\ (\forall)\, n,m \geq n_{p,\varepsilon}.$$
so
$$p(T_n x - Tx) \leq \varepsilon\, p(x),\ (\forall)\, n \geq n_{p,\varepsilon},$$
This implies that
$$\hat{p}(T_n - T) < \varepsilon,\ (\forall)\, n \geq n_{p,\varepsilon}.$$
which prove that $T_n \to T$ in $Q_{\mathcal{P}}(X)$ and $Q_{\mathcal{P}}(X)$ is a sequentially complete m-convex algebra.

Given $(X, \mathcal{P})$, for each $p \in \mathcal{P}$ let $N^p$ denote the null space and $X_p$ the quotient space $X/N^p$. For each $p \in \mathcal{P}$ consider the natural mapping
$$x \to x_p \equiv x + N^p \ (\text{from X to } X_p).$$
It is obvious that $X_p$ is normed space, for each $p \in \mathcal{P}$, with norm defined by
$$\|x_p\|_p = p(x),\ (\forall)\, x \in X.$$
Consider the algebra homomorphism $T \to T^p$ of $Q_{\mathcal{P}}(X)$ into $\mathcal{L}(X_p)$ defined by
$$T^p(x_p) = (Tx)_p,\ (\forall)\, x \in X.$$
This operators are well defined because $T(N^p) \subset N^p$. Moreover, for each $p \in \mathcal{P}$, $\mathcal{L}(X_p)$ is a unital normed algebra and we have
$$\|T_p\|_p = \sup\left\{\|T_p x_p\|_p\ \Big|\ \|x_p\|_p \leq 1 \text{ for } x_p \in X_p\right\} = \sup\{p(Tx)\ |\ p(x) \leq 1 \text{ for } x \in X\}.$$
For $p \in \mathcal{P}$ consider the normed space $(\tilde{X}_p, \|\ \|_p)$ the completition of $(X_p, \|\ \|_p)$. If $T \in Q_{\mathcal{P}}(X)$, then the operator $T^p$ has a unique continuous linear extension $\tilde{T}^p$ on $(\tilde{X}_p, \|\ \|_p)$.

**Definition 1.9.** Let $(X,\mathcal{P})$ a locally convex space and $T \in Q_{\mathcal{P}}(X)$. We say that $\lambda \in \rho(Q_{\mathcal{P}}, T)$ if the inverse of $\lambda I - T$ exists and $(\lambda I - T)^{-1} \in Q_{\mathcal{P}}(X)$.

Spectral sets $\sigma(Q_{\mathcal{P}}, T)$ are defined to be complements of resolvent sets $\rho(Q_{\mathcal{P}}, T)$.

For each $p \in \mathcal{P}$ we denote by $\sigma(X_p, T^p)$ ($\sigma(\tilde{X}_p, \tilde{T}^p)$) the spectral set of the operator $T^p$ in $\mathcal{L}(X_p)$ (respectively the spectral set of $\tilde{T}^p$ in $\mathcal{L}(\tilde{X}_p)$). The resolvent set of the operator $T^p$ in $\mathcal{L}(X_p)$ (respectively the spectral set of $\tilde{T}^p$ in $\mathcal{L}(\tilde{X}_p)$) is denoted by $\rho(X_p, T^p)$ ($\rho(\tilde{X}_p, \tilde{T}^p)$).

**Lemma 1.10.** ([6]) Let $(X,\mathcal{P})$ be a sequentially complete convex space and $T \in Q_{\mathcal{P}}(X)$. Then, T is invertible in $Q_{\mathcal{P}}(X)$ if and only if $\tilde{T}^p$ is invertible in $\mathcal{L}(\tilde{X}_p)$ for all $p \in \mathcal{P}$.

**Corollary 1.11.** ([6]) If $(X,\mathcal{P})$ is a sequentially complete convex space and $T \in Q_{\mathcal{P}}(X)$, then
$$\sigma(Q_{\mathcal{P}}, T) = \cup\{\sigma(X_p, T^p)\ |\ p \in \mathcal{P}\} = \cup\{\sigma(\tilde{X}_p, \tilde{T}^p)\ |\ p \in \mathcal{P}\}.$$

## 2. Spectral radius of quotient bounded operators

Let $(X,\mathcal{P})$ be a locally convex space and $T \in Q_\mathcal{P}(X)$. We said that T is bounded element of the algebra $Q_\mathcal{P}(X)$ if it is bounded element of $Q_\mathcal{P}(X)$ in the sense of G.R.Allan [1]. The class of bounded element of $Q_\mathcal{P}(X)$ is denoted by $(Q_\mathcal{P}(X))_0$.

**Definition 2.1.** If $(X,\mathcal{P})$ is a locally convex space and $T \in Q_\mathcal{P}(X)$ we denote by $r_\mathcal{P}(T)$ the radius of boundness of operator T in $Q_\mathcal{P}(X)$, i.e.

$$r_\mathcal{P}(T) = \inf\{\alpha > 0 \mid \alpha^{-1} T \text{ generates a bounded semigroup in } Q_\mathcal{P}(X)\}.$$

We said that $r_\mathcal{P}(T)$ is the $\mathcal{P}$–spectral radius of the operator T.

Proposition 1.7 (3) implies that for each $\mathcal{P}' \in \mathbf{C}(X)$, $\mathcal{P} \approx \mathcal{P}'$, we have $Q_{\mathcal{P}'}(X) = Q_\mathcal{P}(X)$, so if $\mathcal{H}$ is a Q-equivalence class in $\mathbf{C}(X)$, then

$$r_\mathcal{P}(T) = r_{\mathcal{P}'}(T), \; (\forall) \; \mathcal{P},\mathcal{P}' \in \mathcal{H}.$$

**Proposition 2.2.** ( [1]) If X is a locally convex algebra and $\mathcal{P} \in \mathbf{C}(X)$, then for each $T \in Q_\mathcal{P}(X)$ we have

$$r_\mathcal{P}(T) = \sup\{\limsup_{n\to\infty} (\hat{p}(T^n))^{1/n} \mid p \in \mathcal{P}\}.$$

From real analysis we have following lemma.

**Lemma 2.3.**( [10]) If $(t_n)_n$ is a sequence in $\mathbf{R}^+ \cup \{\infty\}$ then

$$\limsup_{n\to\infty} \sqrt[n]{t_n} = \inf\left\{ v > 0 \;\middle|\; \lim_{n\to\infty} \frac{t_n}{v^n} = 0 \right\}.$$

This lemma implies that for a bounded operator on Banach space we have

$$r(T) = \lim_{n\to\infty} \sqrt[n]{\|T^n\|} = \inf\left\{ v > 0 \;\middle|\; \text{sequence } \left(\frac{T^n}{v^n}\right)_n \text{ converge to zero in operator norm topology} \right\}.$$

If we consider this relation as an alternative definition of the spectral radius, then proposition 2.2 implies that $\mathcal{P}$–spectral radius of a quotient bounded operator can be considered to be natural generalization of the spectral radius of bounded operator on Banach space.

**Proposition 2.4.** If X is a locally convex space and $\mathcal{P} \in \mathbf{C}(X)$, then for each $T \in Q_\mathcal{P}(X)$ we have:
1) $r_\mathcal{P}(T) \geq 0$ and

$$r_\mathcal{P}(\lambda T) = |\lambda| r_\mathcal{P}(T), \; (\forall) \; \lambda \in \mathbf{C},$$

where by convention $0 \infty = \infty$;
2) $r_\mathcal{P}(T) < +\infty$ if and only if $T \in (Q_\mathcal{P}(X))_0$;

3) $r_{\mathcal{P}}(T) = \inf\left\{ \lambda > 0 \mid \lim_{n\to\infty} \dfrac{T^n}{\lambda^n} = 0 \right\}$.

**Proof.** a) From proposition 2.2. results that for each $\lambda \in \mathbf{N}^*$, $\lambda \neq 0$, we have

$$r_{\mathcal{P}}(\lambda T) = \sup\{ \limsup_{n\to\infty} (\hat{p}((\lambda T)^n))^{1/n} \mid p \in \mathcal{P} \} = \sup\{ \limsup_{n\to\infty} (|\lambda| \hat{p}(T^n))^{1/n} \mid p \in \mathcal{P} \} =$$

$$= |\lambda| \sup\{ \limsup_{n\to\infty} (\hat{p}(T^n))^{1/n} \mid p \in \mathcal{P} \} = |\lambda| r_{\mathcal{P}}(T).$$

The case $\lambda = 0$ is obviously.

b) This equivalence is a result of definition of $\mathcal{P}$–spectral radius of the operator T.

c) The equality results directly from proposition 2.2 and lemma 2.3.

**Proposition 2.5.** If X is a locally convex space and $\mathcal{P} \in \mathbf{C}(X)$, then for each $T \in Q_{\mathcal{P}}(X)$ and each $p \in \mathcal{P}$ the sequence $(\hat{p}(T_n))_n$ is convergent and

$$\lim_{n\to\infty} (\hat{p}(T^n))^{1/n} = \inf_{n\geq 1} (\hat{p}(T^n))^{1/n}.$$

**Proof.** If there exists some $n_0 \in \mathbf{N}^*$ such that $\hat{p}(T^{n_0}) = 0$, then

$$\hat{p}(T^n) \leq \hat{p}(T^{n-n_0}) \hat{p}(T^{n_0}) = 0, (\forall) n > 0.$$

Therefore,

$$\lim_{n\to\infty} (\hat{p}(T^n))^{1/n} = \inf_{n\geq 1} (\hat{p}(T^n))^{1/n}.$$

Assume that $\hat{p}(T^n) > 0$ for each $n \in \mathbf{N}^*$. Let $m \in \mathbf{N}^*$ be arbitrary fixed. For each $n \in \mathbf{N}^*$ we consider the relations

$$n = m \cdot q(n) + r(n),$$

where $0 \leq r(n) < m$. Using this notations we have

$$(\hat{p}(T^n))^{1/n} = (\hat{p}(T^{mq(n)+r(n)}))^{1/n} \leq [\hat{p}(T^{mq(n)}) \hat{p}(T^{r(n)})]^{1/n} \leq$$

$$\leq \hat{p}(T^{mq(n)})^{1/n} \hat{p}(T^{r(n)})^{1/n} \leq \hat{p}(T^m)^{q(n)/n} \hat{p}(T^{r(n)})^{1/n},$$

so

$$\limsup_{n\to\infty} (\hat{p}(T^n))^{1/n} \leq (\hat{p}(T^m))^{1/m}.$$

Since $m \in \mathbf{N}^*$ is arbitrary fixed, from previous inequality results that

$$\limsup_{n\to\infty} (\hat{p}(T^n))^{1/n} \leq \inf_{m\geq 1} (\hat{p}(T^m))^{1/m} \leq \liminf_{n\to\infty} (\hat{p}(T^n))^{1/n}.$$

Therefore, the sequence $(\hat{p}(T_n))_n$ is convergent and

$$\lim_{n\to\infty} (\hat{p}(T^n))^{1/n} = \inf_{n\geq 1} (\hat{p}(T^n))^{1/n}.$$

**Corollary 2.6.** If X is a locally convex space and $\mathcal{P} \in \mathbf{C}(X)$, then for each $T \in Q_{\mathcal{P}}(X)$ we have:

a) $r_{\mathcal{P}}(T) = \sup\{ \lim_{n\to\infty} (\hat{p}(T^n))^{1/n} \mid p \in \mathcal{P} \} = \sup\{ \inf_{n\geq 1} (\hat{p}(T^n))^{1/n} \mid p \in \mathcal{P} \}$;

b) $r_{\mathcal{P}}(T) \leq \hat{p}(T), (\forall) p \in \mathcal{P}.$

**Lemma 2.7.** Let $(X,\mathcal{P})$ a locally convex space and $T,S \in Q_{\mathcal{P}}(X)$. If $TS=ST$, then for each $p \in \mathcal{P}$ we have
$$\limsup_{n\to\infty}\left(\hat{p}\left((TS)^n\right)\right)^{1/n} \leq \limsup_{n\to\infty}\left(\hat{p}(T^n)\right)^{1/n} \limsup_{n\to\infty}\left(\hat{p}(S^n)\right)^{1/n}.$$

**Proof.** From inequality
$$\hat{p}(TS) \leq \hat{p}(T)\hat{p}(S),$$
we have
$$\hat{p}\left((TS)^n\right) \leq \hat{p}(T^n)\hat{p}(S^n), \ (\forall)\ n \in \mathbf{N}^*.$$
Therefore,
$$\limsup_{n\to\infty}\left(\hat{p}\left((TS)^n\right)\right)^{1/n} \leq \limsup_{n\to\infty}\left(\hat{p}(T^n)\right)^{1/n} \limsup_{n\to\infty}\left(\hat{p}(S^n)\right)^{1/n}.$$

**Corollary 2.8.** Let $X$ is a locally convex space and $\mathcal{P} \in C(X)$. If $T,S \in Q_{\mathcal{P}}(X)$ such that $TS=ST$, then
$$r_{\mathcal{P}}(TS) \leq r_{\mathcal{P}}(T)\ r_{\mathcal{P}}(S).$$

**Lemma 2.9.** Let $(X,\mathcal{P})$ a locally convex space and $T,S \in Q_{\mathcal{P}}(X)$. If $TS=ST$, then for each $p \in \mathcal{P}$ there exists some seminorm $q \in \mathcal{P}$ such that
$$\limsup_{n\to\infty}\left(\hat{p}\left((T+S)^n\right)\right)^{1/n} \leq \limsup_{n\to\infty}\left(\hat{p}(T^n)\right)^{1/n} + \limsup_{n\to\infty}\left(\hat{p}(S^n)\right)^{1/n}.$$

**Proof.** For each $\varepsilon > 0$ there exists some seminorm index $n_\varepsilon \in \mathbf{N}$ such that
$$\left(\hat{p}(T^m)\right)^{1/m} < \limsup_{n\to\infty}\left(\hat{p}(T^n)\right)^{1/n} + \varepsilon,$$
$$\left(\hat{p}(S^m)\right)^{1/m} < \limsup_{n\to\infty}\left(\hat{p}(S^n)\right)^{1/n} + \varepsilon,$$
for every $m \geq n_\varepsilon$. Therefore, there exists $k \geq 1$ such that
$$\hat{p}(T^m) < k\left(\limsup_{n\to\infty}\left(\hat{p}(T^n)\right)^{1/n} + \varepsilon\right)^m,$$
$$\hat{p}(S^m) < k\left(\limsup_{n\to\infty}\left(\hat{p}(S^n)\right)^{1/n} + \varepsilon\right)^m,$$
for every $m \in \mathbf{N}$. Moreover, $p$ is a submultiplicative seminorm so we have
$$\hat{p}\left((T+S)^m\right) \leq \sum_{k=1}^{m} C_m^k \hat{p}(T^j)\hat{p}(S^{m-j}) \leq$$
$$\leq k^2 \sum_{k=1}^{m} C_m^k \left(\limsup_{n\to\infty}\left(\hat{p}(T^n)\right)^{1/n} + \varepsilon\right)^k \left(\limsup_{n\to\infty}\left(\hat{p}(S^n)\right)^{1/n} + \varepsilon\right)^{m-k} =$$
$$= k^2 \left(\limsup_{n\to\infty}\left(\hat{p}(T^n)\right)^{1/n} + \limsup_{n\to\infty}\left(\hat{p}(S^n)\right)^{1/n} + 2\varepsilon\right)^m.$$
for every $n \geq n_\varepsilon$. Therefore
$$\limsup_{n\to\infty}\left(\hat{p}\left((T+S)^n\right)\right)^{1/n} \leq \limsup_{n\to\infty}\left(\hat{p}(T^n)\right)^{1/n} + \limsup_{n\to\infty}\left(\hat{p}(S^n)\right)^{1/n} + 2\varepsilon$$
Since $\varepsilon > 0$ is arbitrary chosen from previous relation results that

$$\limsup_{n\to\infty}\left(\hat{p}\left((T+S)^n\right)\right)^{1/n} \le \limsup_{n\to\infty}\left(\hat{p}\left(T^n\right)\right)^{1/n} + \limsup_{n\to\infty}\left(\hat{p}\left(S^n\right)\right)^{1/n}.$$

**Corollary 2.10.** Let X is a locally convex space and $\mathcal{P}\in C(X)$. If $T,S\in Q_{\mathcal{P}}(X)$ such that TS=ST, then

$$r_{\mathcal{P}}(T+S) \le r_{\mathcal{P}}(T) + r_{\mathcal{P}}(S).$$

From definition of the $\mathcal{P}$-spectral radius of a quotient bounded operator and previous observation result the following proposition

**Proposition 2.11.** Let X is a locally convex space and $\mathcal{P}\in C(X)$.
1) If $T\in (Q_{\mathcal{P}}(X))_0$, then

$$\lim_{n\to\infty}\frac{T^n}{\lambda^n}=0,\ (\forall)\ |\lambda|>r_{\mathcal{P}}(T);$$

2) If $T\in (Q_{\mathcal{P}}(X))_0$ and $0<|\lambda|<r_{\mathcal{P}}(T)$, then the set $\left\{\dfrac{T^n}{\lambda^n}\right\}_{n\ge 1}$ is unbounded.

3) For each $T\in Q_{\mathcal{P}}(X)$ and every $n>0$ we have

$$r_{\mathcal{P}}(T^n) = r_{\mathcal{P}}(T)^n.$$

**Proposition 2.12.** Let X be a sequentially complete locally convex space and $\mathcal{P}\in C(X)$. If $T\in (Q_{\mathcal{P}}(X))_0$ and $|\lambda|>r_{\mathcal{P}}(T)$, then the Neumann series $\sum_{n=0}^{\infty}\dfrac{T^n}{\lambda^{n+1}}$ converges to $R(\lambda,T)$ (in $Q_{\mathcal{P}}(X)$) and $R(\lambda,T)\in Q_{\mathcal{P}}(X)$.

**Proof.** Let $\beta\in \mathbf{C}$ such that $0<|\beta|<1$ and $r_{\mathcal{P}}(T)<\beta\lambda$. From proposition 2.11(1) we obtain that for each $\varepsilon>0$ and every $p\in\mathcal{P}$, there exists some index $n_{p,\varepsilon}\in \mathbf{N}$, with property

$$\hat{p}\left(\frac{T^n}{(\beta\lambda)^n}\right)<\varepsilon,\ (\forall)\ n\ge n_{p,\varepsilon}.$$

therefore, using corollary 1.4. we obtain

$$p\left(\frac{T^n}{(\beta\lambda)^n}x\right)\le \hat{p}\left(\frac{T^n}{(\beta\lambda)^n}\right)p(x)<\varepsilon\ p(x),\ (\forall)\ n\ge n_{p,\varepsilon},\ (\forall)\ x\in X.$$

Since $0<|\beta|<1$, there exists $n_0\in \mathbf{N}$, such that

$$\sum_{k=n}^{m}|\beta|^k<1,\ (\forall)\ m>n\ge n_0.$$

From previous relations result that for each $\varepsilon>0$ and every $p\in\mathcal{P}$ there exists an index $m_{p,\varepsilon}=\max\{n_{p,\varepsilon},n_0\}\in \mathbf{N}$, for which we have

$$p\left(\sum_{k=n}^{m}\frac{T^k}{\lambda^k}x\right)\le \varepsilon\left(\sum_{k=n}^{m}|\beta|^k\right)p(x)<\varepsilon\ p(x), \qquad (3)$$

for every $m>n\ge m_{p,\varepsilon}$ and every $x\in X$.

Therefore, $\left( \sum_{k=0}^{m} \frac{T^k}{\lambda^{k+1}} x \right)_{m \geq 0}$ is a Cauchy sequence, for each $x \in X$.

But X is sequentially complete, so for every $x \in X$ there exists an unique element $y \in X$ such that

$$y = \lim_{m \to \infty} \sum_{k=0}^{m} \frac{T^k}{\lambda^{k+1}} x.$$

We consider the operator $S: X \to X$ given by

$$S(x) = \lim_{m \to \infty} \sum_{k=0}^{m} \frac{T^k}{\lambda^{k+1}} x, \ (\forall) \ x \in X.$$

It is obvious that S is linear operator. Moreover, from equality

$$\sum_{k=0}^{m} \frac{T^k}{\lambda^{k+1}} (\lambda x - Tx) = x - \frac{T^{m+1}}{\lambda^{m+1}} x, \ (\forall) \ x \in X,$$

result that if $m \to \infty$ then

$$S(\lambda x - Tx) = x, \ (\forall) \ x \in X.$$

so $S(\lambda I - T) = I$. From continuity of the operator T result that

$$ST x = \lim_{m \to \infty} \sum_{k=0}^{m} \frac{T^k}{\lambda^{k+1}} Tx = \lim_{m \to \infty} T\left( \sum_{k=0}^{m} \frac{T^k}{\lambda^{k+1}} x \right) = T\left( \lim_{m \to \infty} \sum_{k=0}^{m} \frac{T^k}{\lambda^{k+1}} x \right) = TS x,$$

for all $x \in X$, therefore

$$S(\lambda I - T) = (\lambda I - T)S = I.$$

The definition of $\mathcal{P}$–spectral radius implies that family $\left( \frac{T^n}{(\beta \lambda)^n} \right)_n$ is bounded in $Q_{\mathcal{P}}(X)$, therefore for every $p \in \mathcal{P}$ there exists a constant $\varepsilon_p > 0$ with property

$$\hat{p}\left( \frac{T^n}{(\beta \lambda)^n} \right) < \varepsilon_p, \ (\forall) \ n \geq 1.$$

Using again corollary 1.4. we have

$$p\left( \frac{T^n}{\lambda^n} x \right) < \varepsilon_p \ |\beta|^n \ p(x), \ (\forall) \ n \geq 1, \ (\forall) \ x \in X,$$

Therefore, for every $p \in \mathcal{P}$ there exists some $\varepsilon_p > 0$ with property

$$p\left( \sum_{k=0}^{m} \frac{T^k}{\lambda^{k+1}} x \right) < \frac{\varepsilon_p}{|\lambda|} \left( \sum_{k=0}^{m} |\beta|^k \right) p(x) < \frac{\varepsilon_p}{|\lambda|} \frac{1}{1-|\beta|} p(x),$$

for every $m \geq 1$ and every $x \in X$, which implies that $S = R(\lambda, T) \in Q_{\mathcal{P}}(X)$.

If we write the relation (3) under form

$$p\left( \sum_{k=0}^{m} \frac{T^k}{\lambda^{k+1}} x - \sum_{k=0}^{n} \frac{T^k}{\lambda^{k+1}} x \right) < \frac{\varepsilon}{|\lambda|} p(x),$$

then for $m \to \infty$ result that for every $\varepsilon > 0$ and every $p \in \mathcal{P}$ there exists some index $n_{p,\varepsilon} \in \mathbf{N}$, such that

$$p\left( Sx - \sum_{k=0}^{n} \frac{T^k}{\lambda^{k+1}} x \right) < \frac{\varepsilon}{|\lambda|} p(x), \ (\forall) \ n \geq n_{p,\varepsilon}, \ (\forall) \ x \in X,$$

This implies that the Neumann series $\sum_{n=0}^{\infty} \frac{T^n}{\lambda^{n+1}}$ converges to $R(\lambda,T)$ in $Q_{\mathcal{P}}(X)$.

**Proposition 2.13.** Let X be a sequentially complete locally convex space and $\mathcal{P} \in \mathbf{C}(X)$. If $T \in Q_{\mathcal{P}}(X)$, then
$$|\sigma(Q_{\mathcal{P}},T)| = r_{\mathcal{P}}(T).$$

**Proof.**
Inequality $|\sigma(Q_{\mathcal{P}},T)| \leq r_{\mathcal{P}}(T)$ is implied by previous proposition.

We prove now the reverse inequality. From corollary 1.11. we have
$$\sigma(Q_{\mathcal{P}},T) = \cup \{\sigma(X_p, T^p) \mid p \in \mathcal{P}\} = \cup \{\sigma(\tilde{X}_p, \tilde{T}^p) \mid p \in \mathcal{P}\}.$$
so, if $|\lambda| > |\sigma(Q_{\mathcal{P}},T)|$, then
$$|\lambda| > |\sigma(\tilde{X}_p, \tilde{T}^p)|, \; (\forall) \, p \in \mathcal{P}.$$
But, $\tilde{X}_p$ is Banach space for each $p \in \mathcal{P}$, therefore
$$|\sigma(\tilde{X}_p, \tilde{T}^p)| = r(\tilde{X}_p, \tilde{T}^p)$$
where $r(\tilde{X}_p, \tilde{T}^p)$ is spectral radius of bounded operator $\tilde{T}^p$ in $\tilde{X}_p$.

This observation implies that for each $p \in \mathcal{P}$ we have $\frac{T^{p\,n}}{\lambda^n} \to 0$ in $\mathcal{L}(\tilde{X}_p)$, which means that for every $p \in \mathcal{P}$ and every $\varepsilon > 0$ there exists $n_{\varepsilon,p} \in \mathbf{N}$, such that
$$\hat{p}\left(\frac{T^n}{\lambda^n}\right) = \left\|\frac{T^{p\,n}}{\lambda^n}\right\|_p < \varepsilon, \; (\forall) \, n \geq n_{\varepsilon,p}.$$
Using proposition 2.4(3) and previous relation we have $r_{\mathcal{P}}(T) \leq |\lambda|$.

But $|\lambda| > |\sigma(Q_{\mathcal{P}},T)|$ is arbitrary chosen, then
$$r_{\mathcal{P}}(T) \leq |\sigma(Q_{\mathcal{P}},T)|.$$

**Definition 2.14.** If X is a locally convex space and $T \in Q_{\mathcal{P}}(X)$, we denote by $\sigma(Q,T)$ the set
$$\cap \{\sigma(Q_{\mathcal{P}},T) \mid \mathcal{P} \in \mathbf{C}(X) \text{ such that } T \in Q_{\mathcal{P}}(X)\}.$$

**Lemma 2.15.** If X is a locally convex space and $T \in Q_{\mathcal{P}}(X)$ then
$$|\sigma(Q,T)| \leq \inf\{r_{\mathcal{P}}(T) \mid \mathcal{P} \in \mathbf{C}(X) \text{ such that } T \in Q_{\mathcal{P}}(X)\}.$$
**Proof.** This is a direct consequence of proposition 2.13.

**Definition 2.16.** An operator T is quotient bounded operator on a locally convex space X if there exists some calibration $\mathcal{P}$ on X such that $T \in Q_{\mathcal{P}}(X)$.

**Observation 2.17.** An operator T is quotient bounded on a locally convex space X if and only if there exists some calibration $\mathcal{P} \in \mathbf{C}(X)$ such that $\hat{p}(T)$ is finit for each $p \in \mathcal{P}$.

**Lemma 2.18.** If X is a locally convex space and T is quotient bounded operator on X, then there exists some principal calibration $\mathcal{P}' \in \mathbf{C}_0(X)$ such that $T \in Q_{\mathcal{P}'}(X)$.

**Proof.** Let $\mathcal{P}$ be a calibration on X such that $T \in Q_{\mathcal{P}}(X)$ and denote by $\mathcal{P}'$ the set of all seminorms given by relations
$$p'(x) = \max_{i=\overline{1,n}} p_i(x), \quad (\forall) \, x \in X,$$
where $p_i \in \mathcal{P}'$, $i = \overline{1,n}$, and $n \in \mathbf{N}$.

Let $p' \in \mathcal{P}'$ be arbitrary chosen. Since $T \in Q_{\mathcal{P}}(X)$, from previous observation and lemma 1.3.3(2) results that
$$p_i(Tx) \le \hat{p}_i(T) \, p_i(x), \quad (\forall) \, x \in X, \, i = \overline{1,n},$$

If $c_{p'} = \max_{i=\overline{1,n}} \hat{p}_i(T)$, then previous inequality implies that
$$p_i(Tx) \le c_{p'} \, p_i(x) \le c_{p'} \, p'(x), \quad (\forall) \, x \in X, \, i = \overline{1,n},$$
so
$$p'(Tx) \le c_{p'} \, p'(x), \quad (\forall) \, x \in X,$$
Therefore, $T \in Q_{\mathcal{P}'}(X)$.

**Lemma 2.19.** If X is a locally convex space and $T \in Q_{\mathcal{P}}(X)$ then
$$\inf \{ r_{\mathcal{P}}(T) \mid \mathcal{P} \in \mathbf{C}_0(X) \text{ such that } T \in Q_{\mathcal{P}}(X) \} =$$
$$= \inf \{ r_{\mathcal{P}}(T) \mid \mathcal{P} \in \mathbf{C}(X) \text{ such that } T \in Q_{\mathcal{P}}(X) \}$$

Proof. Assume that $\mathcal{P} \in \mathbf{C}(X)$ such that $T \in Q_{\mathcal{P}}(X)$.

If $|\lambda| > r_{\mathcal{P}}(T)$, then the family $\left( \dfrac{T^n}{\lambda^n} \right)_{n \ge 0}$ is bounded in $Q_{\mathcal{P}}(X)$, i.e. for every $p \in \mathcal{P}$ there exists $\varepsilon_p > 0$ such that
$$\hat{p} \left( \frac{T^n}{\lambda^n} \right) \le \varepsilon_p, \quad (\forall) \, n \ge 0.$$

Let $\mathcal{P}'$ be the principal calibration associated with the calibration $\mathcal{P}$, i.e. for each $p' \in \mathcal{P}'$ there exists $p_1, ..., p_n \in \mathcal{P}$ such that
$$p' = \max\{p_1, ..., p_n\}.$$
If $\varepsilon_{p'} = \max\{\varepsilon_{p_1}, ..., \varepsilon_{p_n}\}$, then
$$\hat{p}' \left( \frac{T^n}{\lambda^n} \right) \le \varepsilon_{p'}, \quad (\forall) \, n \ge 0.$$
so $|\lambda| > r_{\mathcal{P}'}(T)$. Since $\lambda$ is arbitrary chosen results
$$r_{\mathcal{P}'}(T) \le r_{\mathcal{P}}(T).$$

Therefore,
$$\inf \{ r_{\mathcal{P}}(T) \mid \mathcal{P} \in \mathbf{C}_0(X) \text{ such that } T \in Q_{\mathcal{P}}(X) \} \le$$
$$\le \inf \{ r_{\mathcal{P}}(T) \mid \mathcal{P} \in \mathbf{C}(X) \text{ such that } T \in Q_{\mathcal{P}}(X) \}$$

The reverse inequality is obviously.

**Lemma 2.20.** If X is a locally convex space and $T \in Q_{\mathcal{P}}(X)$ then
$$\sigma(Q,T) = \cap \{ \sigma(Q_{\mathcal{P}},T) \mid \mathcal{P} \in \mathbf{C}_0(X) \text{ such that } T \in Q_{\mathcal{P}}(X) \}.$$
**Proof.** From definition of the set $\sigma(Q,T)$ results the inclusion
$$\sigma(Q,T) \subset \cap \{ \sigma(Q_{\mathcal{P}},T) \mid \mathcal{P} \in \mathbf{C}_0(X) \text{ such that } T \in Q_{\mathcal{P}}(X) \}.$$
If $\lambda \notin \sigma(Q,T)$, such that $\mathcal{P} \in C(X)$ such that $\lambda \in \rho(Q_{\mathcal{P}},T)$, so for every $p \in \mathcal{P}$ we have
$$\hat{p}(R(\lambda,T)) < \infty.$$
Denote by $\mathcal{P}'$ the principal calibration of all seminorms
$$p'(x) = \max_{i=1,n} p_i(x), \ (\forall) \ x \in X,$$
where $p_i \in \mathcal{P}'$, $i = \overline{1,n}$, and $n \in \mathbf{N}$.

Let $p' \in \mathcal{P}'$ be such seminorms. Since $R(\lambda,T) \in Q_{\mathcal{P}}(X)$, the lemma 1.3 (2) implies that
$$p_i(R(\lambda,T)x) \leq \hat{p}_i(R(\lambda,T)) \ p_i(x), \ (\forall) \ x \in X, \ i = \overline{1,n},$$
If $c_{p'} = \max_{i=1,n} \hat{p}_i(T)$, then previous inequality implies that
$$p_i(R(\lambda,T)x) \leq c_{p'} \ p_i(x) \leq c_{p'} \ p'(x), \ (\forall) \ x \in X, \ i = \overline{1,n},$$
so we have
$$p'(R(\lambda,T)x) \leq c_{p'} \ p'(x), \ (\forall) \ x \in X,$$
Therefore, $R(\lambda,T) \in Q_{\mathcal{P}'}(X)$ and $\lambda \notin \sigma(Q_{\mathcal{P}'},T)$, which implies that
$$\cap \{ \sigma(Q_{\mathcal{P}},T) \mid \mathcal{P} \in \mathbf{C}_0(X) \text{ such that } T \in Q_{\mathcal{P}}(X) \} \subset \sigma(Q,T).$$

### 3. Locally bounded operators

**Definition 3.1.** Let $\mathcal{P}$ be a calibration on X. A linear operator $T: X \to X$ is universally bounded on $(X,\mathcal{P})$ if exists a constant $c_0 > 0$ such that
$$p(Tx) \leq c_0 p(x), \ (\forall) \ x \in X.$$
Denote by $B_{\mathcal{P}}(X)$ the collection of all universally bounded operators on $(X,\mathcal{P})$. First we observe that $B_{\mathcal{P}}(X) \subset \mathfrak{L}(X)$.

**Lemma 3.2.** If $\mathcal{P}$ a calibration on X, then $B_{\mathcal{P}}(X)$ is an unital normed algebra with respect to the norm $\| \ \|_{\mathcal{P}}$ defined by
$$\|T\|_{\mathcal{P}} = \sup\{ \hat{p}(T) \mid p \in \mathcal{P} \}, \ (\forall) T \in B_{\mathcal{P}}(X).$$
**Proof.** From definition of $\| \ \|_P$ and lemma 1.7.(2). We have
$$p(Tx) \leq \|T\|_{\mathcal{P}} \ p(x), \ (\forall) x \in X, \ (\forall) \ p \in \mathcal{P}.$$
If $M > 0$ such that
$$p(Tx) \leq M \ p(x), \ (\forall) x \in X, \ (\forall) \ p \in \mathcal{P},$$
hence

$$\hat{p}(T) \le M, \ (\forall) \ p \in \mathcal{P}.$$

So, by definition we have $\|T\|_{\mathcal{P}} \le M$.

**Definition 3.3.** Let X and Y be topological vector spaces. An operator $T: X \to Y$ is said to be locally bounded if it maps some zero neighborhood into a bounded set.

**Observation 3.4.** Let X be a locally bounded spaces.
1) If $T \in L(X)$, then the folowing statements are equivalent:
   a) T is continuous;
   b) T is quotient bounded;
   c) T is locally bounded.
2) If we multiply a locally bounded operator by a continuous operator on the left the product is locally bounded.

The class of local bounded operators is an algebra and it will be usually equipped with the topology of uniform convergence on a zero neighborhood. We say that a sequence $(S_n)$ converge uniformly to zero on a zero neighborhood U if for each zero neighborhood V there exists a positive index $n_0 \in \mathbb{N}$ such that

$$S_n(U) \subset V, \ (\forall) n \ge n_0.$$

In terms of operator seminorms this definition is given by the following definition

**Definition 3.5.** Let X be a locally convex space. We say that a sequence $(S_n)_n$ converges uniformly to zero on some zero-neighborhood if for each principal calibration $\mathcal{P} \in \mathbf{C}_0(X)$ there exists some seminorm $p \in \mathcal{P}$ such that for every $q \in \mathcal{P}$ and every $\varepsilon > 0$ there exists an index $n_{q,\varepsilon} \in \mathbf{N}$, with the property

$$m_{pq}(S_n) < \varepsilon, \ (\forall) \ n \ge n_{q,\varepsilon}$$

A family $\mathcal{G}$ of bounded operators is uniformly bounded on some zero-neighborhood if there exists some seminorm $p \in \mathcal{P}$ such that for every $q \in \mathcal{P}$ there exists $\varepsilon_q > 0$ with the property

$$m_{pq}(S) < \varepsilon_q, \ (\forall) \ S \in \mathcal{G}.$$

**Definition 3.6.** Let T be a locally bounded operator on locally convex space X. We say that $\lambda \in \rho_{lb}(T)$ if there exists a scalar $\alpha$ and a locally bounded operator S on X such that

$$(\lambda I - T)^{-1} = \alpha I + S.$$

The spectral set $\sigma_{lb}(T)$ is defined to be the complement of the resolvent set $\rho_{lb}(T)$.

**Definition 3.7.** Given a linear operator T on a topological vector space X, we consider

$$r_{lb}(T) = \inf \left\{ \nu > 0 \ \middle| \ \frac{T^n}{\nu^n} \to 0 \text{ uniformly on some zero neigborhood} \right\}.$$

**Lemma 3.8.** (lemma 2 [15]) If $T_1$ and $T_2$ are locally bounded operators on X, then there exists a calibration $\mathcal{P}'$ on X such that $T_1, T_2 \in B_{\mathcal{P}'}(X)$.

**Proposition 3.9.** (proposition 8. [3]) Let T be a locally bounded operator on a sequentially complete locally convex space X and $\mathcal{P} \in \mathbf{C}(X)$ such that $T \in Q_{\mathcal{P}}(X)$. If $p \in \mathcal{P}$ such that
$$m_{pq}(T) < \infty, \ (\forall) \ q \in \mathcal{P},$$
then $\rho(X_p, T^p) = \rho(\tilde{X}_p, \tilde{T}^p)$.

**Lemma 3.10.** If X is a locally convex space and $T \in \mathcal{LB}(X)$, then
$$r_{l\,b}(T) = \inf\left\{ \nu > 0 \mid \left(\frac{T^n}{\nu^n}\right)_n \text{ is uniformly bounded on some zero-neighborhood} \right\}.$$

**Proof.** Let $\mathcal{P} \in \mathbf{C}_0(X)$ and
$$r'(T) = \inf\left\{ \nu > 0 \mid \left(\frac{T^n}{\nu^n}\right)_n \text{ is uniformly bounded on some zero-neighborhood} \right\}.$$

Assume that $\nu > r_{l\,b}(T)$, i.e. the sequence $\left(\frac{T^n}{\nu^n}\right)_n$ converges to zero uniformly on some zero neighborhoods. Then, there exists some seminorm $p_1 \in \mathcal{P}$ such that for each $q \in \mathcal{P}$ and every $\varepsilon > 0$ there exists an index $n_{q,\varepsilon} \in \mathbf{N}$, with the property
$$m_{p_1 q}\left(\frac{T^n}{\nu^n}\right) < \varepsilon, \ (\forall) \ n \geq n_{q,\varepsilon}$$

Since T is locally bounded, there exists some seminorm $p_2 \in \mathcal{P}$ such that for each $q \in \mathcal{P}$ we have $m_{p_2 q}(T) < \infty$. But the calibration $\mathcal{P}$ is principal so there exists some seminorm $p_0 \in \mathcal{P}$ such that
$$p_1 \leq p_0 \text{ and } p_2 \leq p_0.$$

Let $\varepsilon > 0$ be arbitrary fixed. Then for each $q \in \mathcal{P}$ we have
$$m_{p_0 q}\left(\frac{T^n}{\nu^n}\right) = \sup\left\{ q\left(\frac{T^n}{\nu^n} x\right) \mid p_0(x) \leq 1 \right\} \leq$$
$$\leq \sup\left\{ q\left(\frac{T^n}{\nu^n} x\right) \mid p_1(x) \leq 1 \right\} = m_{p_1 q}\left(\frac{T^n}{\nu^n}\right), \ (\forall) \ n \geq n_{q,\varepsilon}$$
$$m_{p_0 q}\left(\frac{T}{\nu}\right) = \sup\left\{ q\left(\frac{T}{\nu} x\right) \mid p_0(x) \leq 1 \right\} \leq$$
$$\leq \sup\left\{ q\left(\frac{T}{\nu} x\right) \mid p_2(x) \leq 1 \right\} = m_{p_2 q}\left(\frac{T}{\nu}\right).$$

Moreover, from lemma 1.3(2) results
$$q\left(\frac{T^k}{\nu^k} x\right) \leq m_{p_2 q}\left(\frac{T}{\nu}\right) p_2\left(\frac{T^{k-1}}{\nu^{k-1}} x\right) \leq \ldots \leq m_{p_2 q}\left(\frac{T^k}{\nu^k}\right)^k p_2(x), \ (\forall) \ x \in X,$$
where $k = \overline{1, n_{q,\varepsilon} - 1}$, so by corollary 1.4 we have

$$m_{p_2 q}\left(\frac{T^k}{v^k}\right) \leq m_{p_2 q}\left(\frac{T}{v}\right)^k, \quad k = \overline{1, n_{q,\varepsilon} - 1}.$$

If
$$\alpha_q = \max\left\{\varepsilon, m_{p_2 q}\left(\frac{T}{v}\right), \ldots, m_{p_2 q}\left(\frac{T}{v}\right)^{n_{q,\varepsilon}-1}\right\},$$

then by previous relations results
$$m_{p_0 q}\left(\frac{T^n}{v^n}\right) \leq \alpha_q, \quad (\forall)\, n \in \mathbf{N}.$$

This relation implies that the family $\left(\frac{T^n}{v^n}\right)_n$ is uniformly bounded on some zero-neighborhood, so $r'(T) \leq v$. Since $v > r_{lb}(T)$ is arbitrary chosen we have $r'(T) \leq r_{lb}(T)$.

We prove now the reverse inequality. If $\alpha > r'(T)$, then there exists
$$\beta \in (r'(T), \alpha),$$
such that the sequence $\left(\frac{T^n}{\beta^n}\right)_n$ is uniformly bounded on a zero neighborhood, i.e. the exists some seminorm $p_0 \in \mathcal{P}$ such that for each $q \in \mathcal{P}$ there exists $\beta_q > 0$ for which we have
$$m_{p_0 q}\left(\frac{T^n}{\beta^n}\right) < \beta_q, \quad (\forall)\, n \in \mathbf{N}.$$

Therefore
$$m_{p_0 q}\left(\frac{T^n}{\alpha^n}\right) = \left(\frac{\beta}{\alpha}\right)^n m_{p_0 q}\left(\frac{T^n}{\beta^n}\right) < \left(\frac{\beta}{\alpha}\right)^n \beta_q, \quad (\forall)\, n \in \mathbf{N}.$$

Since $\frac{\beta}{\alpha} < 1$ results that for every $q \in \mathcal{P}$ and each $\varepsilon > 0$ there exists some index $n_{q,\varepsilon} \in \mathbf{N}$ with the property
$$\left(\frac{\beta}{\alpha}\right)^n \beta_q < \varepsilon, \quad (\forall)\, n \geq n_{q,\varepsilon}.$$

Therefore
$$m_{p_0 q}\left(\frac{T^n}{\alpha^n}\right) < \varepsilon, \quad (\forall)\, n \geq n_{q,\varepsilon},$$

so the sequence $\left(\frac{T^n}{\alpha^n}\right)_n$ converges uniformly to zero on some zero neighborhood and $r_{lb}(T) \leq \alpha$.

But $\alpha > r'(T)$ is arbitrary chosen, so $r_{lb}(T) \leq r'(T)$.

**Proposition 3.11.** If X is a locally convex space and $T \in \mathcal{LB}(X)$, then
$$\sigma_{lb}(T) = \sigma(Q, T) = \sigma(T).$$

**Proof.** The inclusion $\sigma(T) \subset \sigma(Q, T)$ is obviously.

If $\lambda \notin \sigma_{lb}(T)$, then there exists $\alpha \in \mathbf{C}$ and $S \in \mathcal{LB}(X)$ such that

$$(\lambda I\text{-}T)^{-1} = \alpha I + S.$$

By lemma 3.8 there exists $\mathcal{P} \in \mathbf{C}(X)$ such that $T, S \in B_{\mathcal{P}}(X)$ and
$$(\lambda I\text{-}T)^{-1} \in B_{\mathcal{P}}(X) \subset Q_{\mathcal{P}}(X),$$
from which results that $\lambda \notin \sigma(Q_{\mathcal{P}}, T)$ and $\sigma(Q,T) \subset \sigma_{lb}(T)$.

We prove the reverse inclusion $\sigma_{lb}(T) \subset \sigma(T)$. There are two different cases.

If X is locally bounded, then by observation 3.4(1) we have $\sigma_{lb}(T) = \sigma(T)$.

Assume that the space X is not locally bounded. If $\lambda \neq 0$, such that $\lambda \notin \sigma(T)$, then the operator $R(\lambda, T)$ is continuous. From equality
$$R(\lambda, T)(\lambda I\text{-}T) = I$$
results
$$R(\lambda, T) = \frac{1}{\lambda} R(\lambda, T) T + \frac{1}{\lambda} I$$

By observation 3.4 (2) the previous equality implies that the operator $\frac{1}{\lambda} R(\lambda, T) T$ is locally bounded, so $\lambda \notin \sigma_{lb}(T)$.

If $\lambda = 0 \notin \sigma(T)$, then the operators $R(0, T) = (-T)^{-1}$ is continuous, so from the equality $I = (-T)^{-1}(-T)$ it results that the identity operator is locally bounded, which contradict the assumption that X is not locally bounded.

Therefore $\lambda = 0$ is in each set considerate ($\sigma(T)$, $\sigma(Q, T)$, respectively $\sigma_{lb}(T)$) in the case that X is not locally bounded.

From previous observations results that $\sigma_{lb}(T) \subset \sigma(T)$.

**Proposition 3.12.** If T is a locally bounded operator on a locally convex space X, then
$$r_{lb}(T) = \inf \{ r_{\mathcal{P}}(T) \mid \mathcal{P} \in \mathbf{C}(X) \text{ such that } T \in Q_{\mathcal{P}}(X) \}.$$

**Proof.** We denote
$$r'(T) = \inf \{ r_{\mathcal{P}}(T) \mid \mathcal{P} \in \mathbf{C}(X) \text{ such that } T \in Q_{\mathcal{P}}(X) \}$$

If $\lambda > r_{lb}(T)$ and $\mathcal{P} \in \mathbf{C}_0(X)$, such that $T \in Q_{\mathcal{P}}(X)$, then there exists $\mu \in (r_{lm}(T), \lambda)$ such that the sequence $\left( \dfrac{T^n}{\mu^n} \right)_{n \geq 1}$ converges to zero on a zero neighborhood, i.e. there exists a seminorm $p \in \mathcal{P}$ such that for each $q \in \mathcal{P}$ and every $\varepsilon > 0$ there exists some index $n_{q,\varepsilon} \in \mathbf{N}$ with the property
$$m_{pq}\left( \frac{T^n}{\mu^n} \right) < \varepsilon, \ (\forall) \ n \geq n_{q,\varepsilon}.$$

Let consider the family of seminorms $\mathcal{Q} = \{ q_m \mid m \geq 1, q \in \mathcal{P} \}$, where
$$q_m(x) = \max\{ m\, p(x), q(x) \}, \ (\forall) \ x \in X.$$

It is obviously that $\mathcal{Q} \in \mathbf{C}(X)$ and T is quotient bounded operator with respect to the calibration $\mathcal{Q}$. If $q_m \in \mathcal{Q}$, then we have
$$q_m\left( \frac{T^n}{\mu^n} x \right) = \max\left\{ m\, p\left( \frac{T^n}{\mu^n} x \right), q\left( \frac{T^n}{\mu^n} x \right) \right\} \leq$$

$$\leq \max\left\{ m\ \hat{p}\left(\frac{T^n}{\mu^n}\right) p(x), m_{pq}\left(\frac{T^n}{\mu^n}\right) p(x)\right\} \leq$$

$$\leq \max\{m\varepsilon\, p(x),\, \varepsilon\, p(x)\} \leq m\varepsilon\, p(x) \leq \varepsilon\, q_m(x),$$

for each $n \geq n_{q_m,\varepsilon} = \max\{n_{p,\varepsilon}, n_{q,\varepsilon}\}$.

Previous relation and corollary 1.4 implies that

$$\hat{q}_m\left(\frac{T^n}{\mu^n}\right) \leq \varepsilon,\ (\forall)\ n \geq n_{q_m,\varepsilon},$$

therefore, since $\frac{\mu}{\lambda} < 1$, results

$$\hat{q}_m\left(\frac{T^n}{\lambda^n}\right) = \left(\frac{\mu}{\lambda}\right)^n \hat{q}_m\left(\frac{T^n}{\mu^n}\right) < \varepsilon,\ (\forall)\ n \geq n_{q_m,\varepsilon}.$$

From proposition 2.4 (3) results that $r_{\mathcal{Q}}(T) \leq \lambda$, so $r'(T) \leq \lambda$. Since $|\lambda| > r_{lm}(T)$ is arbitrary chosen we have $r'(T) \leq r_{lb}(T)$.

Now we prove the opposite inequality. If $\lambda > r'(T)$, then by lemma 2.18 and proposition 2.11 (3) there exists $\mu \in [r'(T), \lambda)$ and $\mathcal{P} \in \mathbf{C}_0(X)$ such that the sequence $\left(\frac{T^n}{\mu^n}\right)_{n\geq 1}$ converges to zero in $Q_{\mathcal{P}}(X)$, i.e. for each $\varepsilon > 0$ and every $p \in \mathcal{P}$ there exists an index $n_{p,\varepsilon} \in \mathbf{N}$ with property

$$\hat{p}\left(\frac{T^n}{\mu^n}\right) < \varepsilon,\ (\forall)\ n \geq n_{p,\varepsilon}.$$

But T is locally bounded, so there exists $q \in \mathcal{P}$ with the property that for each $p \in \mathcal{P}$ there exists $c_p > 0$ such that

$$m_{qp}\left(\frac{T}{\mu}\right) < c_p.$$

Therefore, for each $\varepsilon > 0$ and every $p \in \mathcal{P}$ there exists some index $n_{p,\varepsilon} \in \mathbf{N}$, with the property

$$p\left(\frac{T^{n+1}}{\lambda^{n+1}}x\right) = \left(\frac{\mu}{\lambda}\right)^{n+1} p\left(\frac{T}{\mu}\left(\frac{T^n}{\mu^n}x\right)\right) \leq \left(\frac{\mu}{\lambda}\right)^{n+1} m_{qp}\left(\frac{T}{\mu}\right) q\left(\frac{T^n}{\mu^n}x\right) \leq$$

$$< \left(\frac{\mu}{\lambda}\right)^{n+1} c_p\,\varepsilon\, q(x) < \varepsilon\, q(x),\ (\forall)\ x \in X,$$

for each $n \geq n_{p,\varepsilon}$. From corollary 1.4 results that for each $\varepsilon > 0$ and every $p \in \mathcal{P}$ there exists an index $n_{p,\varepsilon} \in \mathbf{N}$, such that

$$m_{qp}\left(\frac{T^n}{\lambda^n}\right) < \varepsilon,\ (\forall)\ n \geq n_{p,\varepsilon},$$

so the sequence $\left(\frac{T^n}{\lambda^n}\right)_{n\geq 1}$ converges to zero on a zero neighborhood, which means that $r_{lb}(T) \leq \lambda$. Since $|\lambda| > r'(T)$ is arbitrary chosen results that $r_{lb}(T) \leq r'(T)$.

**Lemma 3.13.** Let T be a locally bounded operator on a locally convex space and $\mathcal{P} \in \mathbf{C}(X)$, such that $T \in Q_{\mathcal{P}}(X)$. If $p \in \mathcal{P}$ such that
$$m_{pq}(T) < \infty, \ (\forall) \ q \in \mathcal{P},$$
and $\lambda \in \rho(T)$ have the property
$$\hat{p}(R(\lambda, T)) < \infty,$$
then $\lambda \in \rho(X_p, T^p)$ and
$$R(\lambda, T^p)(x + N_p) = R(\lambda, T)(x) + N^p, \ (\forall) \ x \in X.$$

**Proof.** Denote by $I^p$ the identity operator on the space $X_p$. From definition of the operator $T^p$ results that $\lambda \in \rho(X_p, T^p)$ and is obviously that the operator $S: X^p \to X^p$, given by relation
$$S(x + N_p) = R(\lambda, T)(x) + N_p, \ (\forall) \ x \in X,$$
is linear and continuous on $X_p$. Moreover, from the definition of the operator $T^p$ results the equalities
$$S(\lambda I^p - T^p)(x + N_p) = S((\lambda I - T)^p(x + N_p)) = S((\lambda I - T)(x) + N^p) =$$
$$= R(\lambda, T)(\lambda I - T)(x) + N^p = x + N^p.$$
$$(\lambda I^p - T^p) S(x + N_p) = (\lambda I^p - T^p) S(x + N_p) = (\lambda I^p - T^p)(R(\lambda, T)(x) + N_p) =$$
$$= (\lambda I - T)^p (R(\lambda, T)(x) + N_p) = (\lambda I - T) R(\lambda, T)(x) + N^p = x + N^p,$$
for each $x \in X$. Therefore $S = R(\lambda, T^p)$.

The next theorem was proved by Troitsky [22], but since I consider that this demonstration contains some errors I will give a different demonstration.

**Proposition 3.14.** If T is a locally bounded operator on a sequentially complete locally convex space X, then
$$r_{lb}(T) = |\sigma(T)| = |\sigma_{lm}(T)| = |\sigma(Q, T)|.$$
**Proof.** Since $\sigma_{lm}(T) = \sigma(Q, T) = \sigma(T)$ from the proposition 5.5 [22] results
$$|\sigma(T)| = |\sigma_{lb}(T)| = |\sigma(Q, T)| \leq r_{lb}(T).$$

Let $\mathcal{P} \in \mathbf{C}_0(X)$, arbitrary chosen, such that $T \in Q_{\mathcal{P}}(X)$. Since T is locally bounded there exists some seminorm $p \in \mathcal{P}$ such that
$$m_{pq}(T) < c_q < \infty, \ (\forall) \ q \in \mathcal{P},$$
We prove that
$$\sigma(\tilde{X}_p, \tilde{T}^p) \subset \sigma_{lb}(T) \text{ and } r_{lb}(T) \leq r(\tilde{X}_p, \tilde{T}^p),$$
where $r(\tilde{X}_p, \tilde{T}^p)$ is spectral radius of the operators $\tilde{T}^p$ in algebra $\mathfrak{L}(\tilde{X}_p)$.

We observe that if $p \leq p_1$, then
$$m_{p_1 q}(T) \leq m_{pq}(T) < c_q < \infty, \ (\forall) \ q \in \mathcal{P}.$$

First we will show that $\sigma(\tilde{X}_p, \tilde{T}^p) \subset \sigma_{lb}(T)$. If $\lambda \in \rho_{lb}(T)$, then from proposition 3.11 results $\lambda \in \rho(T)$.

Assume that $r_{lb}(T) < 1$.

I. If $r_{l\,b}(T) < |\lambda|$, then by proposition 5.5 [22] the series $\sum_{n=0}^{\infty} \frac{T^n}{\lambda^{n+1}}$ converges uniformly on zero-neighborhood to $R(\lambda,T)$, i.e. there exists a seminorm $q_0 \in \mathcal{P}$ such that for every seminorm $q \in \mathcal{P}$ and every $\varepsilon > 0$ there exists some index $n_{q,\varepsilon} \in \mathbf{N}$, with property

$$m_{q_0 q}\left( R(\lambda,T) - \sum_{k=0}^{n} \frac{T^k}{\lambda^{k+1}} \right) < \varepsilon, \quad (\forall)\, n \geq n_{q,\varepsilon}.$$

Moreover, from lemma 3.10 the sequence $\left(\frac{T^n}{\lambda^n}\right)_n$ is uniformly bounded on a zero-neighborhood, i.e. there exists some seminorm $q_1 \in \mathcal{P}$ such that for each $q \in \mathcal{P}$ there exists $\beta_q > 0$ for which we have

$$m_{q_1 q}\left(\frac{T^n}{\lambda^n}\right) < \beta_q, \quad (\forall)\, n \in \mathbf{N}.$$

Since the calibration $\mathcal{P}$ is principal results that there exists a seminorm $q_2 \in \mathcal{P}$ such that
$$q_1 \leq q_2 \text{ and } q_0 \leq q_2.$$

From previous relation results that for each $q \in \mathcal{P}$ the folowing inequalities holds

$$q(R(\lambda,T)x) \leq q\left(\left(R(\lambda,T) - \sum_{k=0}^{n_{q\varepsilon}} \frac{T^k}{\lambda^{k+1}}\right)x\right) + q\left(\sum_{k=0}^{n_{q\varepsilon}} \frac{T^k}{\lambda^{k+1}} x\right) \leq$$

$$\leq m_{q_0 q}\left(R(\lambda,T) - \sum_{k=0}^{n_{q\varepsilon}} \frac{T^k}{\lambda^{k+1}}\right) q_0(x) + \sum_{k=0}^{n_{q\varepsilon}} q\left(\frac{T^k}{\lambda^{k+1}} x\right) \leq$$

$$\leq m_{q_0 q}\left(R(\lambda,T) - \sum_{k=0}^{n_{q\varepsilon}} \frac{T^k}{\lambda^{k+1}}\right) q_0(x) + \left(\lambda^{-1} + \sum_{k=1}^{n_{q\varepsilon}} m_{q_1 q}\left(\frac{T^k}{\lambda^{k+1}}\right)\right) q_1(x) \leq$$

$$\leq (\varepsilon + \lambda^{-1} n_{q\varepsilon} \beta_q + \lambda^{-1}) q_2(x), \quad (\forall)\, x \in X.$$

Therefore, for each $q \in \mathcal{P}$ there exists a constant
$$\gamma_q = \varepsilon + \lambda^{-1} n_{q\varepsilon} \beta_q + \lambda^{-1} > 0$$
such that
$$q(R(\lambda,T)x) \leq \gamma_q\, q_2(x), \quad (\forall)\, x \in X.$$

Let $p_1 \in \mathcal{P}$ be a seminorms such that $q_2 \leq p_1$ and $p \leq p_1$. The previous inequality implies that
$$q(R(\lambda,T)x) \leq \gamma_q\, p_1(x), \quad (\forall)\, x \in X,$$
so
$$m_{p_1 q}(R(\lambda,T)) < \infty, \quad (\forall)\, q \in \mathcal{P}, \tag{4}$$
and particularly $\hat{p}(R(\lambda,T)) < \infty$.

Moreover, from the observation we made at the beginning of the demonstration results
$$m_{p_1 q}(T) \leq m_{pq}(T) < c_q < \infty.$$

Without lose the generality of demonstration we can assume that the seminorm $p$ considerate at the beginning of the demonstration is the seminorm $p_1$. Therefore, the conditions of lemma 3.13 are fulfilled, so $\lambda \in \rho(X_p, T^p)$ and by proposition 3.9 we have $\lambda \in \rho(\tilde{X}_p, \tilde{T}^p)$.

If $|\lambda| \leq r_{lb}(T)$, then $r_{lb}(T) < 1 < |\lambda^{-1}|$, and by previous case we have

$$\lambda^{-1} \in \rho(X_p, T^p) = \rho(\tilde{X}_p, \tilde{T}^p).$$

Since $\lambda, \lambda^{-1} \in \rho(T)$ from equalities

$$(\lambda^{-1} - \lambda) R(\lambda, T) (\lambda^{-1} I - T) \left[ (\lambda^{-1} - \lambda)^{-1} I - (\lambda^{-1} I - T)^{-1} \right] =$$
$$= R(\lambda, T) \left[ (\lambda^{-1} I - T) - (\lambda^{-1} - \lambda) I \right] = R(\lambda, T) (\lambda I - T) = I$$

results

$$\left[ (\lambda^{-1} - \lambda)^{-1} I - (\lambda^{-1} I - T)^{-1} \right]^{-1} = (\lambda^{-1} - \lambda) R(\lambda, T) (\lambda^{-1} I - T). \tag{5}$$

which implies that $(\lambda - \lambda^{-1})^{-1} \in \rho(R(\lambda^{-1}, T))$.

The relation (4) ($\lambda^{-1}$ fulfill the condition of the previous case) implies that

$$m_{pq}(R(\lambda^{-1}, T)) < \infty, \ (\forall) \ q \in \mathcal{P},$$

so the lemma 3.13 and the proposition 3.10 implies that the operators

$$(\lambda^{-1} - \lambda)^{-1} \tilde{I}^p - (\lambda^{-1} \tilde{I}^p - \tilde{T}^p)^{-1}$$

is invertible and continuous on the Banach space $\tilde{X}_p$, so

$$\left[ (\lambda^{-1} - \lambda)^{-1} \tilde{I}^p - (\lambda^{-1} \tilde{I}^p - \tilde{T}^p)^{-1} \right]^{-1} \in \mathcal{L}(\tilde{X}_p).$$

Since the relation (5) is equivalently with

$$R(\lambda, T) = (\lambda^{-1} - \lambda)^{-1} \left[ (\lambda^{-1} - \lambda)^{-1} I - (\lambda^{-1} I - T)^{-1} \right]^{-1} (\lambda^{-1} I - T)^{-1}$$

from previous observations results that $\lambda \in \rho(\tilde{X}_p, \tilde{T}^p)$.

II. Assume that $1 \leq r_{lb}(T) < c < \infty$ (T is locally bounded). Then the operator $T_1 = c^{-1} T$ is locally bounded and satisfies the properties specified at the beginning of the demonstration with the respect to the seminorm $p$ and $r_{lb}(T_1) < 1$.

Therefore we are in the conditions of the case I, so

$$c^{-1} \rho_{lb}(T) = \rho_{lb}(T_1) \subset \rho(\tilde{X}_p, \tilde{T}_1^p) = c^{-1} \rho(\tilde{X}_p, \tilde{T}^p),$$

Therefore, we have $\sigma(\tilde{X}_p, \tilde{T}^p) \subset \sigma_{lb}(T)$.

Let show that $r_{lb}(T) \leq r(\tilde{X}_p, \tilde{T}^p)$.

If $\lambda > r(\tilde{X}_p, \tilde{T}^p)$, then there exists $\mu \in (r(\tilde{X}_p, \tilde{T}^p), \lambda)$ such that the sequence $\left( \dfrac{(\tilde{T}^p)^n}{\mu^n} \right)_{n \geq 1}$ converges to zero in $\mathcal{L}(\tilde{X}_p)$, i.e. for every $\varepsilon > 0$ there exists some index $n_{p,\varepsilon} \in \mathbf{N}$ such that

$$\hat{p}\left( \frac{T^n}{\mu^n} \right) = \left\| \frac{(\tilde{T}^p)^n}{\mu^n} \right\|_p < \varepsilon, \ (\forall) \ n \geq n_{p,\varepsilon}.$$

Therefore, for each $q \in \mathcal{P}$ and every $\varepsilon > 0$ we have

$$q\left( \frac{T^{n+1}}{\lambda^{n+1}} x \right) = \left( \frac{\mu}{\lambda} \right)^{n+1} q\left( \frac{T^{n+1}}{\mu^{n+1}} x \right) \leq \left( \frac{\mu}{\lambda} \right)^{n+1} m_{pq}(T) p\left( \frac{T^n}{\mu^{n+1}} x \right) \leq$$

$$\leq \lambda^{-1}\left(\frac{\mu}{\lambda}\right)^n m_{pq}(T) p\left(\frac{T^n}{\mu^n} x\right) \leq \lambda^{-1}\left(\frac{\mu}{\lambda}\right)^n m_{pq}(T) \hat{p}\left(\frac{T^n}{\lambda^n}\right) p(x) \leq$$

$$\leq \lambda^{-1} c_q \left(\frac{\mu}{\lambda}\right)^n \varepsilon\, p(x).$$

But $\frac{\mu}{\lambda} < 1$, so

$$q\left(\frac{T^{n+1}}{\lambda^{n+1}} x\right) < \varepsilon\, p(x)$$

for every $n$ sufficiently large and every $x \in X$.

By corollary 1.4 the previous inequality implies that for every $q \in \mathcal{P}$ and each $\varepsilon > 0$ there exists some index $n_{q,\varepsilon} \in \mathbb{N}$ such that

$$m_{pq}\left(\frac{T^n}{\lambda^n}\right) < \varepsilon,\ (\forall)\, n \geq n_{q,\varepsilon},$$

i.e. the sequence $\left(\frac{T^n}{\lambda^n}\right)_n$ converges uniformly to zero on a zero-neighborhood, so $\lambda > r_{lb}(T)$.

Since $\lambda > r(\tilde{X}_p, \tilde{T}^p)$ is arbitrary chosen results

$$r_{lb}(T) \leq r(\tilde{X}_p, \tilde{T}^p).$$

Therefore, the next statement hold

$$|\sigma(\tilde{X}_p, \tilde{T}^p)| \leq |\sigma(Q, T)| = |\sigma_{lb}(T)| \leq r_{lb}(T) \leq r(\tilde{X}_p, \tilde{T}^p),$$

Since $\tilde{T}^p$ is a bounded operator on the Banach space $\tilde{X}_p$ result

$$|\sigma(\tilde{X}_p, \tilde{T}^p)| = r(\tilde{X}_p, \tilde{T}^p),$$

therefore the proposition is proved.

**Theorem 3.15.** If T is a locally bounded operator on a sequentially complete locally convex space, then the spectral set $\sigma(T)$ is compact

**Proof**.

Let remind to us that $\sigma(T) = \sigma_{lb}(T)$ (proposition 3.11).

If $\lambda \in \rho_{lb}(T)$, then from definition of $\rho_{lb}(T)$ result that there exists some scalar $\alpha \in \mathbb{C}$ and some locally bounded operator S such that

$$(\lambda I - T)^{-1} = \alpha I + S.$$

The lemma 3.8 implies that there exists some calibration $\mathcal{P} \in \mathbf{C}(X)$ such that T, $S \in B_{\mathcal{P}}(X)$. Moreover, we have $(\lambda I - T)^{-1} \in B_{\mathcal{P}}(X)$.

Let $\mu \in \mathbb{C}$ such that $|\mu| < \left\|(\lambda I - T)^{-1}\right\|_{\mathcal{P}}^{-1}$. We will show that $\lambda + \mu \in \rho(T)$.

Since $\left\|\mu(\lambda I - T)^{-1}\right\|_{\mathcal{P}} < 1$, the operatorial series

$$\sum_{k=0}^{\infty} (-\mu)^k (\lambda I - T)^{-(k+1)}$$

converges in operatorial norm of the space $B_{\mathcal{G}}(X)$. By proposition 1 [5] the algebra $B_{\mathcal{G}}(X)$ is complete, so there exists an operator $R(\mu) \in B_{\mathcal{G}}(X)$ such that

$$\sum_{k=0}^{\infty}(-\mu)^k (\lambda I - T)^{-(k+1)} = R(\mu).$$

Using equality

$$((\lambda+\mu)I - T) \ R(\mu) = (\lambda I - T) \ R(\mu) + \mu \ R(\mu) =$$
$$= (\lambda I - T)\left(\sum_{k=0}^{\infty}(-\mu)^k (\lambda I - T)^{-(k+1)}\right) + \mu \left(\sum_{k=0}^{\infty}(-\mu)^k (\lambda I - T)^{-(k+1)}\right) =$$
$$= \sum_{k=0}^{\infty}(-\mu)^k (\lambda I - T)^{-k} - \sum_{k=0}^{\infty}(-\mu)^{k+1} (\lambda I - T)^{-(k+1)} = I$$

$$R(\mu)((\lambda+\mu)I - T) = \mu R(\mu) + R(\mu)(\lambda I - T) =$$
$$= \mu \left(\sum_{k=0}^{\infty}(-\mu)^k (\lambda I - T)^{-(k+1)}\right) + \left(\sum_{k=0}^{\infty}(-\mu)^k (\lambda I - T)^{-(k+1)}\right)(\lambda I - T) =$$
$$= -\sum_{k=0}^{\infty}(-\mu)^{k+1} (\lambda I - T)^{-(k+1)} + \sum_{k=0}^{\infty}(-\mu)^k (\lambda I - T)^{-k} = I$$

we have

$$\lambda + \mu \in \rho(B_{\mathcal{G}}, T) \subset \rho(T).$$

Therefore

$$\{\beta \mid |\beta - \lambda| < \|(\lambda I - T)^{-1}\|_{\mathcal{G}}^{-1}\} \subset \rho(T),$$

Since $\lambda \in \rho_{lm}(T)$ is an arbitrary chosen result that the set $\rho(T)$ is open, so the spectral set $\sigma(T)$ is closed. The spectral set $\sigma(B_{\mathcal{G}}, T)$ is compact, so from the inclusion $\sigma(T) \subset \sigma(B_{\mathcal{G}}, T)$ results that the set $\sigma(T)$ is compact.

**Observation 3.16.** The function

$$\mu \to R(\mu) = ((\lambda+\mu)I - T)^{-1},$$

is analytic at the point $\mu = 0$.

**Corollary 3.17.** Let T be a locally bounded operator on the space X. If $\lambda \in \rho(T)$ and $d(\lambda)$ is the distance from $\lambda$ to the spectrum $\sigma(T)$, then for each calibration $\mathcal{G} \in \mathbf{C}(X)$, with the property $(\lambda I - T)^{-1}$, $T \in B_{\mathcal{G}}(X)$, we have $\|(\lambda I - T)^{-1}\|_{\mathcal{G}} \geq \dfrac{1}{d(\lambda)}$.

**Proof.**
If $\mathcal{G} \in \mathbf{C}(X)$ is a calibration such that $(\lambda I - T)^{-1}$, $T \in B_{\mathcal{G}}(X)$, then in the demonstration of previous theorem we showed that for each scalar $\mu \in \mathbf{C}$ with the property $|\mu| < \|(\lambda I - T)^{-1}\|_{\mathcal{G}}^{-1}$ we have

$$\lambda + \mu \in \rho(B_{\mathcal{G}}, T) \subset \rho(T),$$

so $\|(\lambda I - T)^{-1}\|_{\mathcal{G}}^{-1} \leq d(\lambda)$.